\documentclass{amsart}

\usepackage{amscd}
\usepackage[T1]{fontenc}
\usepackage[latin1]{inputenc}
\usepackage[pdfauthor={Stefan Kebekus and Thomas Peternell}, 
pdftitle={A refinement of Stein factorization and deformations of surjective morphisms},
pdfkeywords={algebraic variety, surjective morphism, deformation of morphism, Hom-scheme, Stein factorization}, 
pdfview={Fit}]{hyperref}
\usepackage[arrow,curve,matrix]{xy}
\usepackage{epsfig}
\newcommand{\sE}{\mathcal{E}}
\newcommand{\sF}{\mathcal{F}}
\newcommand{\sG}{\mathcal{G}}

\renewcommand{\O}{\mathcal{O}}

\DeclareMathOperator{\Aut}{Aut}
\DeclareMathOperator{\Chow}{Chow}
\DeclareMathOperator{\codim}{codim}
\DeclareMathOperator{\Hilb}{Hilb}
\DeclareMathOperator{\Hom}{Hom}
\DeclareMathOperator{\Image}{Image}

\DeclareMathOperator{\Pic}{Pic}
\DeclareMathOperator{\red}{red}

\DeclareMathOperator{\Sing}{Sing}

\DeclareMathOperator{\vertical}{vert}

\newcommand{\ilabel}[1]{\newcounter{#1}\setcounter{#1}{\value{enumi}}}
\newcommand{\iref}[1]{\setcounter{enumi}{\value{#1}}\labelenumi}

\theoremstyle{plain}    
\newtheorem{thm}{Theorem}[section]
\newtheorem{defn}[thm]{Definition}

\numberwithin{equation}{thm}
\numberwithin{figure}{section}
\theoremstyle{plain}    
\newtheorem{cor}[thm]{Corollary}
\newtheorem{lem}[thm]{Lemma}

\theoremstyle{plain}    
\newtheorem{prop}[thm]{Proposition}
\newtheorem{proclaim-special}[thm]{\specialthmname}

\theoremstyle{remark}
\newtheorem{fact}[thm]{Fact}
\newtheorem{rem}[thm]{Remark}
\theoremstyle{remark}    
\newtheorem{claim}[equation]{Claim} 
\newtheorem*{claim*}{Claim} 
\newtheorem{notation}[thm]{Notation} 
\newtheorem{assumption}[thm]{Assumption}
\newtheorem{question}[thm]{Question}

  \def\factor#1.#2.{\left. \raise 2pt\hbox{$#1$} \right/\hskip -2pt\raise -2pt\hbox{$#2$}}
  \renewcommand{\labelenumi}{(\thethm.\arabic{enumi})}

\begin{document}

\title[Refinement of Stein factorization and deformations of
morphisms]{A refinement of Stein factorization and deformations of
  surjective morphisms}

\date{\today}

\author{Stefan Kebekus}
\address{Stefan Kebekus, Mathematisches Institut, Universit\"at zu
  K\"oln, Weyertal 86--90, 50931 K\"oln, Germany}
\email{\href{mailto:stefan.kebekus@math.uni-koeln.de}{stefan.kebekus@math.uni-koeln.de}}
\urladdr{\href{http://www.mi.uni-koeln.de/~kebekus}{http://www.mi.uni-koeln.de/$\sim$kebekus}}

\author{Thomas Peternell}
\address{Thomas Peternell, Institut f\"ur Mathematik, Universit\"at
  Bayreuth, 95440~Bayreuth, Germany}
\email{\href{mailto:thomas.peternell@uni-bayreuth.de}{thomas.peternell@uni-bayreuth.de}}
\urladdr{\href{http://btm8x5.mat.uni-bayreuth.de/mathe1}{http://btm8x5.mat.uni-bayreuth.de/mathe1}}

\thanks{Both authors were supported in part by the
  Forschungsschwerpunkt ``Globale Methoden in der komplexen Analysis''
  of the Deutsche Forschungsgemeinschaft. A part of this paper was
  worked out while Stefan Kebekus visited the Korea Institute for
  Advanced Study.  He would like to thank Jun-Muk Hwang for the
  invitation.}

\begin{abstract}
  This paper is concerned with a refinement of the Stein
  factorization, and with applications to the study of deformations of
  morphisms. We show that every surjective morphism $f : X \to Y$
  between normal projective varieties factors canonically via a finite
  cover of $Y$ that is étale in codimension one.  This ``maximally
  étale factorization'' satisfies a strong functorial property.
  
  It turns out that the maximally étale factorization is stable under
  deformations, and naturally decomposes an étale cover of the
  Hom-scheme into a torus and into deformations that are relative with
  respect to the rationally connected quotient of the target $Y$. In
  particular, we show that all deformations of $f$ respect the
  rationally connected quotient of $Y$.
\end{abstract}

\maketitle

 \setcounter{tocdepth}{1}
 \tableofcontents

\section{Introduction and statement of results}

Throughout this paper, we consider surjective morphisms between
algebraic varieties and their deformations. To fix notation, we use
the following assumption.

\begin{assumption}\label{ass:basic}
  $f: X \to Y$ will always denote a surjective holomorphic map between normal
  complex-projective varieties.
\end{assumption}

The main method that we introduce is a refinement of the Stein
factorization: we show that $f$ factors canonically via a finite cover
of $Y$ that is étale in codimension one.  This ``maximally étale
factorization'' satisfies a strong functorial property which is
defined in Section~\ref{sec:intromaxetale} below and turns out to be
stable under deformations of $f$.

We employ the maximally étale factorization for a study of the
deformation space $\Hom(X,Y)$ and show that an étale cover of the
$\Hom$-scheme naturally decomposes into a torus and into deformations
that are relative with respect to the maximally rationally connected
fibration of the target $Y$. In particular, we show that all
deformations of $f$ respect the rationally connected quotient of $Y$.

These result are summarized and properly formulated below.

\subsection{The maximally étale factorization}
\label{sec:intromaxetale}

Under the Assumptions~\ref{ass:basic}, suppose that there exists a
factorization $f$,
\begin{equation}
  \label{eq:fact_of_f}
  \xymatrix{ {X} \ar[r]_{\alpha} \ar@/^0.3cm/[rr]^{f} & {Z}
    \ar[r]_{\beta} & {Y}}
\end{equation}
where $\beta$ is finite and étale in codimension 1, i.e. étale outside
a set of codimension $\geq 2$.

\begin{defn}\label{def:maxetale}
  We say that a factorization $f = \beta \circ \alpha$ as
  in~\eqref{eq:fact_of_f} is \emph{maximally étale} if the following
  universal property holds: for any factorization $f = \beta' \circ
  \alpha'$, where $\beta': Z' \to Y$ is finite and étale in
  codimension 1, there exists a morphism $\gamma: Z \to Z'$ such that
  such that the following diagram commutes:
  $$
  \xymatrix{ {X} \ar[r]_{\alpha} \ar@/^0.3cm/[rr]^{f} \ar@{=}[d] &
    {Z} \ar[d]_{\gamma} \ar[r]_{\beta} & {Y}  \ar@{=}[d] \\
    {X} \ar[r]^{\alpha'} \ar@/_0.3cm/[rr]_{f} &
    {Z'}\ar[r]^{\beta'} & {Y}} 
  $$
\end{defn}

\begin{rem}
  It follows immediately from the definition that a maximally étale
  factorization of a given morphism $f$ is unique up to isomorphism if
  it exists. Theorem~\ref{thm:char} describes the uniqueness in more
  detail.
\end{rem}

The existence of the maximally étale factorization is established by
the following theorem, which we prove in Section~\ref{sec:pfthm0}.

\begin{thm}\label{thm:main0}
  Let $f : X \to Y$ be a surjective morphism between normal projective
  varieties. Then there exists a maximally étale factorization.
\end{thm}

We will later in Section~\ref{sec:char} describe the maximally étale
factorization in terms of the positivity of the push-forward sheaf
$f_*(\O_X)$.

\begin{rem}
  We have already remarked that the maximally
  étale factorization yields a natural refinement of the Stein
  factorization. More precisely, we can say that a surjection $f: X
  \to Y$ of normal projective varieties decomposes as follows.
  $$
  \xymatrix{ {X} \ar[rr]_{\text{conn.~fibers}}
    \ar@/^0.3cm/[rrrrr]^{f} & & {W} \ar[r]_{\text{finite}} & {Z}
    \ar[rr]_{\text{max.~étale}} & & {Y}}
  $$
\end{rem}

\subsection{Stability of the factorization under deformations}

Let $f=\alpha\circ \beta$, as in Diagram~\eqref{eq:fact_of_f} denote
the Stein factorization. If $f' : X \to Y$ is any deformation of $f$,
it is a classical fact that $f'$ again factors via $\beta$ ---see
Section~\ref{sec:stabilityofstein} for brief review.  We will show
that a similar, and somewhat stronger, property also holds for the
maximally étale factorization.  To formulate this stability result
precisely, we introduce the following notation.

\begin{notation}
  If $g: A \to B$ is any morphism between projective varieties, let
  $\Hom_{g}(A, B) \subset \Hom(A,B)$ be the connected component that
  contains $g$.
\end{notation}

The stability result is then formulated as follows.

\begin{thm}\label{thm:main0a}
  In the setup of Theorem~\ref{thm:main0}, let $f = \beta \circ \alpha
  $ be the maximally étale factorization as in
  Diagram~\eqref{eq:fact_of_f}.  Then the natural
  morphism between the reduced $\Hom$-spaces,
  $$
  \begin{array}{rccc}
    \eta : & \Hom_{\alpha}(X, Z)_{\red} & \to & \Hom_{f}(X, Y)_{\red} \\
    & \alpha' & \mapsto & \beta \circ \alpha'
  \end{array}
  $$
  is proper and surjective. The induced morphism $\tilde \eta$
  between the normalizations is étale.  If $f' \in \Hom_f(X,Y)$ is any
  deformation of $f$, then $f'$ factors via $Z$, and has $Z$ as
  maximally étale factorization.
\end{thm}

We prove Theorem~\ref{thm:main0a} in Section~\ref{sec:pfthm0a}.

\subsection{Decomposition of the $\Hom$-scheme}

We recall the main result of \cite{HKP03}, where deformations of
morphisms with non-uniruled targets were studied. Using the language
of Section~\ref{sec:intromaxetale}, this is formulated as follows.

\begin{thm}[\protect{\cite[thm.~1.2]{HKP03}}]\label{thm:HKP}
  Under the Assumptions~\ref{ass:basic}, suppose that $Y$ is not
  covered by rational curves. If 
  $$
  \xymatrix{
    X  \ar[r]_{\alpha} \ar@/^0.3cm/[rr]^{f} &
    Z  \ar[r]_{\beta} &
    Y
  }
  $$
  denotes the maximally étale factorization, and if $\Aut^0(Z)$ is
  the maximal connected subgroup of the automorphism group of $Z$,
  then $\Aut^0(Z)$ is an Abelian variety, and the natural morphism
  $$
  \factor \Aut^0(Z) . \Aut(Z/Y) \cap \Aut^0(Z) . \to \Hom_f(X,Y)
  $$
  is an isomorphism of schemes.
  
  In particular, all deformations of surjective morphisms $X \to Y$
  are unobstructed, and the associated components of $\Hom_f(X,Y)$ are
  smooth Abelian varieties. \qed
\end{thm}

At the other extreme, if $Y$ is rationally connected, partial
descriptions of the $\Hom$-scheme are known ---the results of
\cite[thm.~1]{HM03} and \cite[thm.~3]{HM04} assert that whenever $Y$
is a Fano manifolds of Picard-number 1 whose variety of minimal
rational tangents is finite, or not linear, then all deformations of
$f$ come from automorphisms of $Y$. This covers all examples of Fano
manifolds of Picard-number one that one encounters in practice.

If $Y$ is covered by rational curves, but not rationally connected, we
consider the maximally rationally connected fibration $q_Y : Y
\dasharrow Q_Y$ which is explained in more detail in
Section~\ref{sec:MRC}. Using the maximally étale factorization, we
will show that an étale cover of $\widetilde \Hom_f(X,Y)$, the
normalization of the space $\Hom_f(X,Y)$, can be decomposed into a
torus and a space of deformations that are relative with respect
$q_Y$. We recall the notion of a relative deformation first.

\begin{notation}\label{not:reldef1}
  We call the subvariety
  $$
  H^f_{\vertical} := \{ f' \in \Hom_f(X,Y)_{\red} \,|\, q_Y \circ f' = q_Y \circ f \}
  $$
  the ``space of relative deformations of $f$ over $q_Y$''.
\end{notation}

The following theorem will then be shown in Section~\ref{sec:decomp}.

\begin{thm}\label{thm:main}
  Under the Assumption~\ref{ass:basic}, let
  $$
  \xymatrix{ X \ar[r]_{\alpha} \ar@/^0.3cm/[rr]^{f} & Z
    \ar[r]_{\beta} & Y }
  $$
  be the maximally étale factorization, and $T \subset \Aut^0(Z)$ a
  maximal compact Abelian subgroup. Then the following holds:
  \begin{enumerate}
  \item\label{thm:m1} There exists a normal variety $\tilde H$ and an
    étale morphism
    $$
    T \times \tilde H \to \widetilde \Hom_f(X,Y)
    $$
    that maps $\{e\} \times \tilde H$ to the preimage of
    $H^f_{\vertical}$.
    
  \item\label{thm:m2} If $Y$ is smooth or if $f$ is itself maximally
    étale, then $\{e\} \times \tilde H$ surjects onto the preimage of
    $H^f_{\vertical}$.
  \end{enumerate}
\end{thm}

In the setup of Theorem~\ref{thm:main}, it need not be true that
$\Aut^0(Z)$ is itself an Abelian variety. The existence of a maximal
compact Abelian subgroup $T \subset \Aut^0(Z)$ is briefly discussed in
Fact~\ref{fact:maxabelian} on page~\pageref{fact:maxabelian} below.

\begin{rem}
  The assertion of Theorem~\ref{thm:main} is weaker than
  Theorem~\ref{thm:HKP} in the sense that it does not make any
  statement about the scheme-structure of $\Hom_f(X,Y)$. The reason is
  that the maximal rationally connected fibration $q_Y$ need not be a morphism,
  and that there is no good deformation space for rational maps
  between fixed varieties.  Theorem~\ref{thm:main} can certainly be
  straightened if one assumes additionally that $q_Y$ is regular.
\end{rem}

\subsection{Acknowledgments}

The authors would like to thank Ivo Radloff and Eckhard Viehweg for a
number of discussions. It was in these discussion that the notion of
``maximally étale'' evolved.

\section{Known Facts}

The proofs of our main results rely on a number of facts scattered
throughout the literature. For the reader's convenience, we have
gathered these here.

\subsection{Stability of Stein factorization under deformation}
\label{sec:stabilityofstein}

Consider the Stein factorization of $f$,
\begin{equation}
  \label{diag:stein}
  \xymatrix{ {X} \ar[rr]_{g \text{: conn. fibers}} \ar@/^0.3cm/[rrrr]^{f} & & {W_0}
    \ar[rr]_{h \text{: finite}} & & {Y}}.
\end{equation}
We will later need to know that any deformation of $f$ still has $g: X
\to W_0$ as Stein factorization. While this is probably
well-understood, we were unable to find a good reference for the
universal properties of Stein factorization, and include a full proof.

\begin{prop}[Stability of Stein factorization under deformation]\label{prop:stability-of-stein}
  The canonical composition morphism
  $$
  \begin{array}{rccc}
    \nu: & \Hom_h(W_0, Y) & \to & \Hom_f(X,Y) \\
    & h' & \mapsto & h' \circ g
  \end{array}
  $$
  is bijective. In particular, the morphism between the normalized
  $\Hom$-schemes is isomorphic.
\end{prop}

The proof of Proposition~\ref{prop:stability-of-stein}, which we give
below, makes use of the following two lemmas, which assert that a
deformation of a morphism with connected fibers does not change the
fibers, and that a surjective morphism between normal spaces is
determined up to isomorphism by its set-theoretical fibers.

\begin{lem}[Invariance of fibers under deformation]\label{lem:invarianceoffibers}
  Let $T$ be a smooth curve and $(f_t)_{t \in T} : X \to Y$ be a
  family of surjective morphisms between projective varieties. Assume
  that for all $t \in T$, the map $f_t$ has connected fibers. Then the
  set-theoretical fibers of $f_t$ are independent of $t$. More
  precisely, for all $x \in X$ and all $s,t \in T$, we have
  $$
  f_t^{-1}f_t(x) = f_s^{-1}f_s(x).
  $$
\end{lem}
\begin{proof}
  Choose an ample line bundle $L \in \Pic(X)$. Observe that two points
  $x, y \in X$ are contained in the same $f_t$-fiber if and only if there exists
  a curve $C \subset X$ that contains both $x$ and $y$, and satisfies
  $c_1 (f_t^*(L)).C=0$. Likewise, $x$ and $y$ are in the same
  $f_s$-fiber if and only if there is a curve $C \subset X$ with $x, y \in C$ and
  $c_1 (f_s^*(L)).C=0$.
  
  On the other hand, since $f_t$ and $f_s$ are homotopic, the Chern
  classes of the pull-back bundles agree,
  $$
  c_1 (f_t^*(L)) = c_1(f_s^*(L)).
  $$
  This shows the claim.
\end{proof}

\begin{lem}\label{lem:samefibers}
  Let $f_1: X \to W_1$ and $f_2: X \to W_2$ be two surjective
  morphisms between normal spaces. If for all $x \in X$ the
  set-theoretical fibers of $f_1$ and $f_2$ agree, i.e.~ if
  $f_1^{-1}(f_1(x)) = f_2^{-1}(f_2(x))$, then there exists a
  commutative diagram as follows.
  $$
  \xymatrix{
    & X \ar@/_.2cm/[dl]_{f_1} \ar@/^.2cm/[dr]^{f_2} \\
    W_1 \ar@{<->}[rr]_{\phi \text{: isomorphic}} && W_2 
  }
  $$
\end{lem}
\begin{proof}
  The morphisms $f_1$, $f_2$ give rise to a morphism from $\iota: X
  \to W_1 \times W_2$, $\iota(x) = (f_1(x),f_2(x))$, and we obtain a
  commutative diagram as follows.
  $$
  \xymatrix{
    & X \ar@/_.2cm/[dl]_{f_1} \ar@/^.2cm/[dr]^{f_2} \ar[d]^{\iota} \\
    W_1 & W_1 \times W_2 \ar[r]_{p_2} \ar[l]^{p_1} & W_2 }
  $$
  The assumptions that the fibers of $f_1$ and $f_2$ agree implies
  that the restrictions $p_i|_{\iota (X)}$ of the morphisms $p_1$ and
  $p_2$ to the image $\iota(X)$ are bijective. Since we are working
  over $\mathbb C$, Zariski's main theorem then implies that the
  restrictions $p_1|_{\iota(X)}$ and $p_2|_{\iota(X)}$ are even
  isomorphic. We can therefore view $\iota(X)$ as the graph of an
  invertible morphism $\phi := p_2|_{\iota(X)} \circ
  (p_1|_{\iota(X)})^{-1}$, which yields the claim.
\end{proof}

\begin{proof}[Proof of Proposition~\ref{prop:stability-of-stein}]
  The injectivity of $\nu$ is obvious because $g$ is surjective. Since
  $\Hom_f(X,Y)$ is connected, to prove surjectivity, it suffices to
  show that any morphism $\gamma_f : T \to \Hom_f(X,Y)$ from a smooth
  irreducible curve $T$ can be lifted to a curve $\gamma_h: T \to
  \Hom_h(W_0, Y)$ such that $\gamma_f = \nu \circ \gamma_h$.
  
  To this end, let
  $$
  \begin{array}{rccc}
    F : & X \times T & \to & Y\times T \\
    & (x,t) & \mapsto & (f_t(x), t)
  \end{array}
  $$
  be the proper product morphism of the universal map and the
  identity, and consider the Stein factorization
  \begin{equation}
    \label{eq:stf}
  \xymatrix{ {X \times T} \ar[rr]_{G  \text{: conn. fibers}} \ar@/^0.5cm/[rrrr]^{F} & & {W}
    \ar[rr]_{H  \text{: finite}} & & {Y \times T}}.
  \end{equation}
  Lemma~\ref{lem:invarianceoffibers} on the invariance of fibers
  implies that the morphisms $g \times Id_T : X \times T \to W_0
  \times T$ and $G$ have the same fibers.  Lemma~\ref{lem:samefibers}
  then asserts that there exists an isomorphism $\phi$ such that the
  factorization~\eqref{eq:stf} extends to a commutative diagram as
  follows.
  \begin{equation}
    \label{eq:stein}
    \xymatrix{ {X\times T}  \ar@/_0.3cm/[rrd]_{g \times Id_T} \ar[rr]_{G \text{: conn. fibers}}
      \ar@/^0.5cm/[rrrr]^{F} & & {W} \ar[rr]_{H \text{:
          finite}}  &
      & {Y \times T}  \ar[r]^{p_Y} & Y\\
      & & {W_0 \times T} \ar[u]^{\phi}
    }
  \end{equation}
  We use Diagram~\eqref{eq:stein} to define the morphism $\gamma_h: T
  \to \Hom_h(W_0, Y)$ by
  $$
  \begin{array}{rccc}
    \gamma_h(t) : & W_0 & \to & Y \\
    & w & \mapsto & p_Y(H(\phi(w,t)))
  \end{array}
  $$
  It follows then from the commutativity of
  Diagram~\eqref{eq:stein} that $\gamma_f = \nu \circ \gamma_h$. This
  proves Proposition~\ref{prop:stability-of-stein}.
\end{proof}

\subsection{The rational quotient}
\label{sec:MRC}

Recall from \cite[chap.~IV]{K96} or \cite[sec.~4]{Debarre01} that an
irreducible projective variety $X$ is \emph{rationally connected} if
any two sufficiently general points can be joined by a {\it single} rational
curve.  Moreover $X$ is \emph{rationally chain connected} if two
general points can be joined by a connected chain of rational curves.

\begin{rem} 
  If $X$ is smooth, then $X$ is rationally connected if and only if
  $X$ is rationally chain connected \cite{K96}. If $X$ is singular,
  this need no longer be true.  For instance, if $X$ is the cone over
  an elliptic curve, then $X$ is of course rationally chain connected,
  but not rationally connected.
\end{rem}

One of the most important features of uniruled varieties is the
existence of a rationally connected quotient, introduced by Campana
and Koll\'ar-Miyaoka-Mori.

\begin{defn}
  Let $V$ be a normal variety and $r_V : V \dasharrow R_V$ a dominant
  rational map to a normal variety. The map $r_V$ is called a
  \emph{maximal rationally chain connected fibration}, if for all very
  general points $v \in V$, the closure of the fiber through $v$,
  $$
  R(v) := \overline{r_V^{-1}(r_V(v))},
  $$
  is the largest rationally chain connected subvariety of $V$ that
  contains $v$.
\end{defn}

The existence of a maximal rationally chain connected fibration is
established by Campana (even in the Kähler case) and
Kollár-Miyaoka-Mori, see \cite{K96} and \cite{Debarre01}. Campana uses
the notation ``rational quotient''.

\begin{fact}\label{fact:mrcf}
  Let $V$ be a normal projective variety. Then there exists a
 maximal rationally chain connected fibration $r_V : V \dasharrow R_V$, with the
  additional property that the quotient map $r_V$
  is \emph{almost holomorphic}, i.e.~there exists a dense open subset
  $V^0 \subset V$ such that the restriction $r_V|_{V^0}$ is a proper
  morphism.
\end{fact}

Notice that Kollár-Miyaoka-Mori \cite{K96} and Debarre
\cite{Debarre01} already put the property to be almost holomorphic
into the definition of a maximal rationally chain connected fibration.
We include this into the next notion.

We will need to consider a different variant of a rational quotient
coming from the fact that for singular varieties rational
connectedness and rational chain connectedness do not coincide.

\begin{defn} 
  Let $V$ be a normal variety and $q_V : V \dasharrow Q_V$ a dominant
  rational map to a normal variety. The map $q_V$ is called a
  \emph{maximal rationally connected fibration}, if for all very general points
  $v \in V$, the closure of the fiber through $v$,
  $$
  R(v) := \overline{q_V^{-1}(q_V(v))},
  $$
  is the largest rationally connected subvariety of $V$ that
  contains $v$.
\end{defn}

\begin{prop}
  If $V$ is a normal projective variety, then there exists a maximal
  rationally connected fibration $q_V: V \dasharrow Q_V$.
\end{prop}
\begin{proof}
  Let $\eta: \tilde V \to V$ be a desingularization, and $r_{\tilde
    V}: \tilde V \dasharrow R_{\tilde V}$ a maximal rationally chain
  connected fibration. Then set $Q_V = R_{\tilde V}$ and $q_V :=
  r_{\tilde V} \circ \eta^{-1}$.
\end{proof}

Notice that there is a factorization $Q_V \dasharrow R_V.$ Of course
both fibration are unique up to birational equivalence, so that we
speak of ``the'' maximal rationally (chain) connected fibration.

\begin{rem}
  If $V$ is singular, a maximal rationally connected fibration of $V$
  is \emph{not} necessarily the maximal rationally chain fibration.
  E.g., if $X$ is the cone over an elliptic curve, then the maximal
  rationally connected fibration maps to the elliptic curve, whereas
  the maximal rationally chain connected fibration maps to a point.
  Further,the maximal rationally connected fibration cannot \emph{not}
  necessarily be taken to be almost holomorphic.
\end{rem}
  
It is a crucial fact shown by Graber, Harris and Starr that the base
of a maximal rationally chain connected fibration, hence also of a
maximal rationally connected fibration is itself not uniruled.

\begin{fact}[\protect{\cite[cor.~1.4]{GHS03}}]
  If $q_V : V \dasharrow Q_V$ is a maximal rationally (chain)
  connected fibration, then $Q_V$ is not uniruled.
\end{fact}

The maximal rationally (chain) connected fibration described in the
literature is determined only up to birational equivalence. It is,
however, easy to see that there is a canonical choice.

\begin{prop}\label{prop:canonrc}
  Let $V$ be a normal projective variety. Then there exists a
  canonical maximal rationally (chain) connected fibration $q_V: V
  \dasharrow Q_V$, with the following property: the automorphism group
  $\Aut(V)$ stabilizes the indeterminacy locus of $q_V$, and has a
  natural action on $Q_V$ such that $q_V$ is equivariant wherever it
  is defined.
\end{prop}
\begin{proof}
  Let $q : V \dasharrow Q$ be any maximal rationally (chain)
  fibration. The universal property of the cycle space than yields a
  rational map as follows:
  $$
  \begin{array}{rccc}
    q'_V : & V & \dasharrow & \Chow(V) \\
    & v & \mapsto & R(v).
  \end{array}
  $$
  
  This construction has two important features. For one, observe that
  the morphism $q'_V$ is independent of the particular choice of the
  rationally connected quotient $q$.  Secondly, if $x\in V$ is a very
  general point, and $g \in \Aut(V)$ is any automorphism, then
  $g(R(x))$ is again rationally (chain) connected.  In particular, we
  have that $g(R(x))) = R(g(x))$.  This already shows that the natural
  action of $\Aut(V)$ on $\Chow(V)$ stabilizes the image of $q'_V$ and
  makes $q'_V$ equivariant.  The proof is thus finished if we let
  $Q_V$ be the normalization of the closure of the image of $q'_V$,
  and $q_V : V \dasharrow Q_V$ be the lifting that comes from the
  universal property of the normalization.
\end{proof}

\begin{notation}\label{not:alwayscanaon}
  For the rest of this paper, if we discuss ``the'' maximal rationally
  (chain) connected fibration of a variety, we always mean the canonic
  construction given in Proposition~\ref{prop:canonrc}.
\end{notation}

We will later need to consider subsheaves of the tangent sheaf $T_V$
that are relative over the rationally connected quotient wherever this
is well-defined.

\begin{defn}\label{def:vertical}
  Let $V$ be a normal projective variety, and let $q_V : V \dasharrow
  Q_V$ be the rationally connected quotient. Further, suppose that $C$
  is a normal variety and $\iota : C \to V$ a morphism whose image is
  not contained in the singular locus of $V$, and not contained in the
  indeterminacy locus of $q_V$. If $\sF \subset \iota^*(T_V)$ is a
  reflexive subsheaf of the pull-back of the tangent sheaf, we say
  that $\sF$ is \emph{vertical with respect to the rationally
    connected quotient}, if $\sF$ is contained in $\iota^*(T_{V|Q_V})$
  at the general point of $C$.
  
  Likewise, a morphism of reflexive sheaves $\sF \to \iota^*(T_V)$ is
  \emph{vertical with respect to the rationally connected quotient} if
  the double dual of its image is. An infinitesimal deformation of
  $\iota$, i.e.~an element $\sigma \in \Hom( \iota^*(\Omega^1_V), \O_C
  )$, is vertical if the restriction of $\sigma$ to general points of
  $C$ corresponds to a section in $\iota^*(T_{V|Q_V})$.
\end{defn}

\subsection{General curves in projective varieties}

We will later need to consider the Harder-Narasimhan filtration of the
tangent sheaf $T_X$. By Mehta-Ramanathan's theorem, it suffices to
discuss the filtration of the restriction to a general complete
intersection curve, whose definition we recall now.

\begin{defn}\label{def:MR}
  If $X$ is normal, we consider \emph{general complete intersection
    curves} in the sense of Mehta-Ramanathan, $C \subset X$. These are
  reduced, irreducible curves of the form $C = H_1 \cap \cdots \cap
  H_{\dim X-1}$, where the $H_i \in |m_i \cdot H|$ are general, the
  $L_i \in \Pic(X)$ are ample and the $m_i \in \mathbb N$ large
  enough, so that the Harder-Narasimhan-Filtration of $T_X$ commutes
  with restriction to $C$.  If the $L_i$ are chosen, we also call $C$
  a \emph{general complete intersection curve with respect to $(L_1,
    \ldots, L_{\dim X-1})$.}  
\end{defn}
  
We refer to \cite{Flenner84} and \cite{Langer04} for a discussion and
an explicit bound for the $m_i$ that appear in
Definition~\ref{def:MR}.

If $X$ is a normal variety, $q_X : X \dasharrow Q_X$ the maximal rationally
connected fibration and $C \subset X$ a general complete intersection
curve, then $C$ intersects neither the singular locus of $X$, nor the
indeterminacy locus of $q_X$. It makes therefore sense ask if a
subsheaf $\sF_C \subset T_X|_C$ is vertical with respect to the
rationally connected quotient. The following important criterion is a
refinement of Miyaoka's characterization of uniruled varieties. It
appeared first implicitly in \cite[9.0.3]{secondasterisque}, but see
\cite[rem.~4.8]{KLT05}.

\begin{fact}[\protect{\cite[cor.~1.4]{KLT05}}]\label{fact:HN3}
  Let $X$ be a normal projective variety with maximal rationally
  connected fibration $q_X: X \dasharrow Q_X$. If $C \subset X$ is a
  general complete intersection curve and $\sF_C \subset T_X \vert C$
  a locally free and ample subsheaf, then $\sF_C \subset T_{X/{Q_X}}
  \vert C$. \qed
\end{fact}

\subsection{Finite morphisms}

Let $f: X \to Y$ be a surjective, finite morphism between normal
varieties. The push-forward of the structure sheaf $f_*(\O_X)$ is then
a torsion free sheaf on $X$, which is locally free where $f$ is flat,
i.e.~away from the singularities of $X$ and $Y$. Much of our
argumentation is based on an analysis of the positivity properties of
$f_* (\O_X)$.

\begin{notation}
  Let $X_{\Sing}$ and $Y_{\Sing}$ denote the singular loci, and set
  \begin{align*}
    X^0 & = X \setminus (X_{\Sing} \cup f^{-1}(Y_{\Sing})) \\
    Y^0 & = Y \setminus (Y_{\Sing} \cup f(X_{\Sing})) = f(X^0).
  \end{align*}
  Then $\codim X\setminus X^0 = \codim Y\setminus Y^0 \geq 2$.
\end{notation}

\begin{fact}\label{fact:splitting}
  The trace morphism $tr: f_* (\O_{X^0}) \to \O_{Y^0}$ gives rise to a
  splitting
  $$
  f_* (\O_{X^0}) \cong \O_{Y^0} \oplus \sE^{\vee}_f 
  $$
  where $\sE^{\vee}_f$ is a locally free sheaf on $Y^0$. Let
  $\sE_f$ be the dual\footnote{The use of the 'dual' follows
    historical conventions.  We use it to be consistent with the
    literature we cite.} of $\sE^{\vee}_f$. \qed
\end{fact}

The following result on the positivity of $\sE_f$ appeared only
recently. We have however learned from E.~Viehweg that it is
implicitly contained in much older works of Fujita.

\begin{fact}[\protect{\cite[Thm.~A of the appendix by 
    Lazarsfeld]{PS00}}]\label{fact:branchFormula} Let $C \subset Y^0$
  be a complete curve that is not contained in the branch locus of
  $f$. Then $\sE_f|_C$ is a nef vector bundle on $C$.  It has degree 0
  if and only if $f$ is unbranched along $C$. \qed
\end{fact}

\begin{cor}\label{cor:etalitecriterion}
  If $C \subset Y$ is a general complete intersection curve, then
  $f_*(\O_X)|_C$ is of degree 0 if and only if $f$ is étale in
  codimension 1.  \qed
\end{cor}

As a consequence of the projection formula, $f_* f^* (\sF) = f_*(\O_X)
\otimes \sF$, we obtain that if $\sF$ is any coherent sheaf on $Y$,
then there is a natural direct sum decomposition
\begin{equation}\label{eq:split}
  H^0\bigl(X^0, f^*(\sF)\bigr) \cong H^0\bigl(Y^0, \sF\bigr) \oplus \Hom_{Y^0}\bigl(\sE_f,\sF|_{Y^0}\bigr)  
\end{equation}

\begin{notation}\label{not:decomp}
  In the setup of this section, if $\sigma \in H^0(X, f^*(\sF))$, let
  $\sigma = \sigma'_f + \sigma''_f$ be the decomposition that is
  associated with the splitting~\eqref{eq:split}.
\end{notation}

\section{Existence of a max.~étale factorization, Proof of Theorem~\ref{thm:main0}}
\label{sec:pfthm0}

We will in this section prove the existence of a maximally étale
factorization for surjective morphisms between normal projective
varieties. Since the proof is somewhat long, we subdivide it into a
number of steps. We maintain the notation and the assumptions of
Theorem~\ref{thm:main0} throughout.

The strategy of proof follows~\cite{HKP03}: we construct the
factorization using a suitable subsheaf of $f_*(\O_X)$.

\subsection{Reduction to the case of a finite morphism}

Using the Stein factorization of the morphism $f$, we can assume
without loss of generality that $f$ is actually finite.

\subsection{The Harder-Narasimhan-Filtration}

Choose an ample line bundle $H \in \Pic(Y)$, and consider the
associated Harder-Narasimhan-Filtration of $f_*(\O_X)$,
\begin{equation}
  \label{eq:HNF}
  0 = \sF_0 \subset \sF_1 \subset \cdots \subset \sF_{r-1} \subset
  \sF_r = f_*(\O_X)
\end{equation}

\begin{lem}\label{lem:sF1}
  The degree of $\sF_1$ with respect to $H$ is zero, 
  $$
  \deg_H (\sF_1) := c_1(\sF_1)\cdot c_1(H)^{\dim Y-1} = 0.
  $$
  If $\sG \subset f_*(\O_X)$ is any coherent subsheaf with
  $\deg_H(\sG) = 0$, then $\sG \subset \sF_1$.
\end{lem}
\begin{proof}
  Consider the splitting $f_*(\O_X) = \O_Y \oplus \sE^*$. Since a
  general complete intersection curve $C \subset Y$ is not contained
  in the branch locus of $f$, Lazarsfeld's result,
  Fact~\ref{fact:branchFormula}, asserts that $\sE^*|_C$ is an
  anti-nef vector bundle. This in turn implies that no subsheaf of
  $f_*(\O_X)$ has positive degree. Since $\O_Y \subset f_*(\O_X)$ is a
  subsheaf of $\deg_H (\O_Y)= 0$, either $f_*(\O_X)$ is $H$-semistable
  and $\deg_H f_*(\O_X) = 0$, or the degree of the maximally
  destabilizing subsheaf $\sF_1$ is zero. The first statement thus
  follows.
  
  The second statement is void if $f_*(\O_X)$ is semistable (so that
  $\sF_1 = f_*(\O_X)$). We can thus assume that $f_*(\O_X)$ is not
  semistable, and that we are given a coherent subsheaf $\sG \subset
  f_*(\O_X)$ with $\deg_H(\sG) = 0$. Consider the image of $\sF_1$ and
  $\sG$ under the addition map,
  $$
  + : \sF_1 \oplus \sG \to f_*(\O_X).
  $$
  The image sheaf again has non-negative degree and must therefore
  be contained in the maximally destabilizing subsheaf $\sF_1$. This
  proves that $\sG \subset \sF_1$.
\end{proof}

\begin{lem}
  The $\O_Y$-algebra structure on $f_*(\O_X)$ induces on $\sF_1$ the
  structure of a sheaf of $\O_Y$-subalgebras.
\end{lem}
\begin{proof}
  Since $\sF_1$ is a sheaf of $\O_Y$-modules which contains $\O_Y$, it
  solely remains to verify that $\sF_1$ is closed under the
  multiplication map
  $$
  m : f_*(\O_X) \otimes f_*(\O_X) \to f_*(\O_X).
  $$
  In other words, we need to check that the associated map
  $$
  m' : \sF_1 \otimes \sF_1 \to \factor f_*(\O_X).\sF_1.
  $$
  is constantly zero. Again, if $\sF_1 = f_*(\O_Y)$, there is
  nothing to show. Otherwise, observe that $\sF_1 \otimes \sF_1$ is
  semistable with slope $\mu(\sF_1 \otimes \sF_1) = 0$ so that
  $\factor f_*(\O_X).\sF_1.$ contains a subsheaf $\sG$ with $\deg_H
  \sG = 0$. By Lemma~\ref{lem:sF1} this subsheaf must vanish, hence
  $m' = 0.$
\end{proof}

\subsection{Construction of the factorization, end of proof}

Since $\sF_1$ is a coherent sheaf of $\O_Y$-algebras, the morphism $f$
now automatically factorizes via $Z := {\rm \bf Specan}(\sF_1)$.
\begin{equation}
  \label{eq:maxetale}
  \xymatrix{ {X} \ar[r]_{\alpha} \ar@/^0.3cm/[rr]^{f} & {Z}
    \ar[r]_{\beta} & {Y}}
\end{equation}
Since $\beta$ is proper and affine, it is clear that it must be
finite. We will now show that $Z$ is normal, that $\beta$ is étale in
codimension 1, and that it is indeed maximally étale.

\begin{lem}
  The intermediate variety $Z$ is normal.
\end{lem}
\begin{proof}
  Let $\eta: \tilde Z \to Z$ be the normalization morphism. The
  universal property of normalization then yields a further
  factorization
  $$
  \xymatrix{ {X} \ar[r]_{\tilde \alpha} \ar@/^0.3cm/[rrr]^{f} &
    {\tilde Z} \ar[r]_{\eta} & {Z} \ar[r]_{\beta} & {Y}}.
  $$
  Accordingly, we obtain a sequences of subsheaves of
  $\O_Y$-algebras,
  $$
  \sF_1 = \beta_*(\O_Z) \subset (\beta \circ \eta)_*(\O_{\tilde Z})
  \subset f_*(\O_X).
  $$
  Since $\eta$ is isomorphic away from a proper subset, the
  quotient $\mathcal Q := \factor (\beta \circ \eta)_*(\O_{\tilde
    Z}).\sF_1.$ either vanishes, or is a torsion sheaf. But since
  $\sF_1$ is saturated in $f_*(\O_X)$, the quotient $\mathcal Q$
  cannot be non-zero torsion. This shows that $(\beta \circ
  \eta)_*(\O_{\tilde Z}) = \sF_1$ and therefore $Z = \tilde Z$.
\end{proof}

\begin{lem}
  The morphism $\beta$ is étale in codimension 1.
\end{lem}
\begin{proof}
  By Corollary~\ref{cor:etalitecriterion}, to prove the assertion, it
  is equivalent to show $\deg (\beta_*(\O_Z)|_C) = 0$.  Since
  $\beta_*(\O_Z) = \sF_1$, and since $\deg (\beta_*(\O_Z)|_C) = \deg_H
  (\sF_1)$, this follows from the first statement of
  Lemma~\ref{lem:sF1}.
\end{proof}

\begin{lem}
  The factorization~\eqref{eq:maxetale} is maximally étale.
\end{lem}
\begin{proof}
  Let $f = \beta' \circ \alpha'$ be any factorization via an
  intermediate variety $Z'$ which is étale in codimension 1 over $Y$.
  The push-forward $\sG := \beta'_*(\O_{Z'}) \subset f_*(\O_X)$ is
  then a subsheaf of $\O_Y$-algebras. If $C \subset Y$ is a general
  complete intersection curve associated with the polarization $H$, then
  $\sG$ is locally free along $C$, and Fact~\ref{fact:branchFormula}
  asserts that $\deg (\sG|_C) = 0$. In other words, we have that
  $\deg_H(\sG) = 0$, and the second statement of Lemma~\ref{lem:sF1}
  implies that $\sG \subset \sF_1 = \beta_*(\O_Z)$.
\end{proof}

This ends the proof of Theorem~\ref{thm:main0}. \qed

\begin{rem}
  If $X$ and $Y$ are smooth, the maximally étale factorization
  \begin{equation}
    \xymatrix{ {X} \ar[r]_{\alpha} \ar@/^0.3cm/[rr]^{f} & {Z}
      \ar[r]_{\beta} & {Y}}.
  \end{equation}
  can more easily be constructed as follows: The subgroup
  $f_*(\pi_1(X)) \subset \pi_1(Y)$ has finite index, and therefore
  determines a finite étale cover $g: \tilde Y \to Y$ such that $f$
  factors via $g$. As the map $\pi_1(X) \to \pi_1(\tilde Y)$ must
  necessarily be onto, the factorization via $g$ is maximal.
\end{rem}

\section{Characterization of the maximally étale factorization}
\label{sec:char}

We will later need to characterize the maximally étale factorization
among all factorizations in terms of positivity properties of the
push-forward sheaf $\beta_*(\O_Z)$. The construction of the maximally
étale factorization in the previous section almost immediately yields
the following.

\begin{thm}\label{thm:char}
  Let $f : X \to Y$ be a surjection between normal projective
  varieties with maximally étale factorization
  \begin{equation}
    \label{eq:5maxet}
    \xymatrix{ {X} \ar[r]_{\alpha} \ar@/^0.3cm/[rr]^{f} & {Z}
      \ar[r]_{\beta} & {Y}}.
  \end{equation}
  Let $H \in \Pic(Y)$ be an arbitrary polarization and $C \subset Y$
  an associated general complete intersection curve. Then
  \begin{enumerate}
  \item \ilabel{I511} the push-forward $\beta_*(\O_Z)$ is the
    maximally destabilizing subsheaf of $f_*(\O_X)$ with respect to
    the polarization $H$, and $Z = {\rm \bf Specan}(\beta_*(\O_Z))$,
    and
    
  \item \ilabel{I512} if we set $\mathcal Q := \factor
    f_*(\O_X).\beta_*(\O_Z).$, then $\mathcal Q^\vee|_C$ is an ample
    vector bundle on $C$.
  \end{enumerate}
\end{thm}
\begin{proof}
  Statement~\iref{I511} is a direct corollary to the proof of
  Theorem~\ref{thm:main0}. In fact, in Section~\ref{sec:pfthm0}, we
  have chosen one particular polarization $H \in \Pic(Y)$, and
  constructed $Z$ as the {\bf Specan} of the maximally destabilizing
  subsheaf $\sF_1 \subset f_*(\O_X)$. While $\sF_1$ could a priori
  depend on the choice of $H$, the universal property of the maximally
  étale factorization shows that it actually does not: if $\beta' : Z'
  \to Y$ is another maximally étale factorization, constructed with
  respect to another polarization $H'$, the universal property of $Z$
  yields the inclusion $\beta'_*(\O_{Z'}) \subset \beta_*(\O_Z) =
  \sF_1$.  Analogously, we obtain that $\sF_1 \subset
  \beta'_*(\O_{Z'})$. This shows statement~\iref{I511}.
  
  It follows from Fact~\ref{fact:branchFormula} that
  $f_*(\O_X)^\vee|_C$ is nef. Since $\beta_*(\O_Z)|_C$ has degree 0,
  it is a standard fact that $\mathcal Q^\vee|_C$ is nef ---see
  \cite[prop.~1.2(8)]{CP91}.  On the other side, Lemma~\ref{lem:sF1}
  implies that $\mathcal Q|_C$ has no subsheaf of semi-positive
  degree. As a consequence, its dual $\mathcal Q^\vee|_C$ has no
  quotient of semi-negative degree.  Hartshorne's characterization
  \cite{Hartshorne71} of ample vectorbundles then implies that
  $\mathcal Q^\vee|_C$ is ample, as claimed.
\end{proof}

\begin{cor}\label{cor:EAntiAmple}
  In the setup of Theorem~\ref{thm:char}, if $H' \in \Pic(Z)$ is any
  polarization, and $C' \subset Z$ an associated general complete
  intersection curve, then the dual of the restriction $\left. \factor
    \alpha_*(\O_X).\O_Z.\right|_{C'}$ is an ample vector bundle on
  $C'$.
\end{cor}
\begin{proof}
  It follows from the universal property of the maximally étale
  factorization~\eqref{eq:5maxet} that the maximally étale
  factorization of $\alpha : X \to Z$ is the identity on $Z$,
  $$
  \xymatrix{ {X} \ar[r]_{\alpha} \ar@/^0.3cm/[rr]^{\alpha} & {Z}
    \ar[r]_{Id} & {Z}}.
  $$
  The claim then follows from Theorem~\ref{thm:char}\iref{I512}.
\end{proof}

\begin{question}
  The Harder-Narasimhan filtration~\eqref{eq:HNF} of $f_*(\O_X)$ that
  is discussed on page~\pageref{eq:HNF} obviously depends on the
  choice of the line bundle $H$. As we have seen in
  Theorem~\ref{thm:char}, it turns out \emph{a posteriori} that the
  maximally destabilizing subsheaf $\sF_1$ does \emph{not} depend on
  $H$.  Are there \emph{a priori} arguments to see that in our setup
  the maximally destabilizing subsheaf is independent of the
  polarization?
\end{question}

\section{Stability under deformations, Proof of Theorem~\ref{thm:main0a}}
\label{sec:pfthm0a}

Throughout the present section we maintain the notation and the
assumptions of Theorem~\ref{thm:main0a}. Again we subdivide the
lengthy proof into steps: after a reduction to the case where $f$ is
finite, we prove the surjectivity of the composition morphism $\eta$
and the étalité of its lift to the normalizations separately.

\subsection{Reduction to the case of a finite morphism}

As an immediate consequence of the stability of Stein factorization
under deformation, Proposition~\ref{prop:stability-of-stein}, we can
replace $X$ with its Stein factorization. We will therefore assume
without loss of generality for the remainder of the present
Section~\ref{sec:pfthm0a} that $f$ is finite.

\subsection{Properness and surjectivity of the composition morphism $\eta$}

The proof of surjectivity is technically a little awkward because the
connected spaces $\Hom_f(X,Y)$ and $\Hom_{\alpha}(X,Z)$ need not be
irreducible. Thus, as a first step, we show that for any irreducible
component of $H \subset \Hom_f(X,Y)$ and any $\alpha'$ with
$\eta(\alpha') \in H$, the component $H$ is the proper image of a
suitable component in $\Hom_{\alpha}(X,Z)$ that contains $\alpha'$.
Surjectivity and properness are then deduced in
Corollaries~\ref{cor:prop} and \ref{cor:surjprop} below.

\begin{prop}\label{prop:p11}
  Let $f' = \alpha' \circ \beta \in \Hom_f(X,Y)_{\red}$ be any
  morphism that factors via $\beta$. Further, let $H_{f'} \subset
  \Hom_f(X,Y)_{\red}$ be an irreducible component that contains $f'$.
  Then there exists a component $H_{\alpha'} \subset \Hom_{\alpha}(X,
  Z)_{\red}$ that contains $\alpha'$ such that $\eta(H_{\alpha'}) =
  H_{f'}$ and such that the restriction $\eta|_{H_{\alpha'}}$ is
  proper.
\end{prop}

\begin{proof}
  Let $\tilde H_{f'}$ be the universal cover of a desingularization of
  $H_{f'}$, and let $\tilde f' \in \tilde H_{f'}$ be a point that maps
  to $f'$. Using that $f'$ factors via $\beta$, we obtain the
  following fibered product diagram:
  $$
  \xymatrix{
    & F \ar[r] \ar[d]_{\tilde \beta} & Z \ar[d]^{\txt{\scriptsize $\beta$ \\\scriptsize étale in codim.~1}} \\
    \{\tilde f'\}\times X \ar[r] \ar@/^.3cm/[ur]  & \tilde H_{f'} \times X \ar[r]^(.6){\mu} \ar[d]_{\txt{\scriptsize $p_2$ \\\scriptsize projection}} & Y  \\
    & X }
  $$
  
  \begin{claim}\label{claim:0absc}
    The morphism $\tilde \beta$ is also étale in codimension 1.
  \end{claim}
  
  \begin{proof}[Proof of Claim~\ref{claim:0absc}]
    Let $R \subset Y$ be the minimal closed set $R$ such that $\beta$
    is étale away from $R$.  Since étale morphisms are stable under
    base change, we only need to show that $\tilde R := \mu^{-1}(R)$
    is of codimension $\geq 2$ in $\tilde H_{f'} \times X$. This will
    be done by showing that for all $\tilde g \in \tilde H$, the
    intersection $\tilde R \cap (\{ \tilde g\} \times X)$ is of
    codimension $\geq 2$ in $\{ \tilde g\} \times X$.
    
    To this end, let $g \in H$ be the image of $\tilde g$. If we
    identify $\{ \tilde g\} \times X$ with $X$ in the obvious way, it
    is then clear that
    $$
    \tilde R \cap (\{ \tilde g\} \times X) = g^{-1}(R).
    $$
    Since $g$ is a deformation of the finite, surjective morphism
    $f$, $g$ is likewise finite and surjective, and
    Claim~\ref{claim:0absc} follows.
  \end{proof}
  
  As a next step in the proof of Proposition~\ref{prop:p11}, let $F^0$
  be the normalization of the irreducible component that contains the
  image of $\{\tilde f'\}\times X$, and let $\tilde \beta^0 : F^0 \to
  \tilde H_{f'} \times X$ be the obvious restriction.

  \begin{claim}\label{claim:iso}
    The morphism $\tilde \beta^0$ is biholomorphic.
  \end{claim}

  \begin{proof}[Proof of Claim~\ref{claim:iso}]
    If $x \in X$ is a general point, set
    $$
    \tilde H_x := p_2^{-1}(x) \text{\quad and \quad} F_x^0 :=
    (\tilde \beta^0)^{-1}(\tilde H_x) \cap F^0.
    $$
    By Seidenberg's theorem \cite{Manaresi82}, $F_x^0$ is normal,
    and the existence of the section $\{\tilde f'\}\times X \cong X
    \to F^0$ implies that $F_x^0$ is irreducible.
    
    Since $\tilde H_{f'} \times X$ is normal, Claim~\ref{claim:0absc}
    now asserts that $\tilde \beta^0$ is étale in codimension 1. Since
    $x$ is general, this is also true for the restriction
    $$
    \tilde \beta^0|_{F^0_x} : F^0_x \to \tilde H_x.
    $$
    But because $\tilde H_x$ is smooth, Zariski-Nagata's theorem on
    the purity of the branch locus, \cite[thm.~3.1]{SGA1}, implies
    that $\tilde \beta|_{F_x}$ must in fact be étale. Since $\tilde
    H_x$ is simply connected, it must be isomorphic. Consequence: the
    finite morphism $\tilde \beta^0$ is bimeromorphic. By the analytic
    version of Zariski's main theorem, \cite[prop.~14.7]{Remmert94},
    $\tilde \beta^0$ is isomorphic. This shows Claim~\ref{claim:iso}.
  \end{proof}
  
  To end the proof of Proposition~\ref{prop:p11}, observe that
  Claim~\ref{claim:iso} shows the existence of a morphism $F^0 \cong
  \tilde H_{f'} \times X \to Z$. The universal property of the
  $\Hom$-scheme thus yields a morphism $\nu: \tilde H_{f'} \to
  \Hom(X,Z)_{\red}$. It follows immediately from the construction that
  $\nu(f')=\alpha'$. Better still, we obtain a diagram
  \begin{equation}\label{eq:factviaeta}
    \xymatrix{
      \tilde H_{f'} \ar@{->>}[d]_{\txt{\scriptsize desing. and \\\scriptsize univ. cover}} \ar[r]^(.3){\nu} & \Hom(X,Z)_{\red} \ar@/^0.3cm/[dl]^{\ \eta} \\
      H_{f'} }
  \end{equation}
  This shows that there exists a component $H_{\alpha'}$ which
  contains $\alpha'$ and surjects onto $H_{f'}$. The properness of
  $\eta|_{H_{\alpha'}}$ follows from Diagram~\eqref{eq:factviaeta}
  because $\eta$ is quasi-finite. This ends the proof of
  Proposition~\ref{prop:p11}.
\end{proof}

\begin{cor}\label{cor:prop}
  Let $f' = \alpha' \circ \beta \in \Hom_f(X,Y)_{\red}$ be any
  morphism that factors via $\beta$. Further, let $H_{\alpha'} \subset
  \Hom_{\alpha'}(X,Z)_{\red}$ be any irreducible component that contains
  $\alpha'$. Then there exists a component $H_{f'} \subset \Hom_{f}(X,
  Y)_{\red}$ that contains $f'$ such that $\eta(H_{\alpha'}) = H_{f'}$
  and such that the restriction $\eta|_{H_{\alpha'}}$ is proper.
\end{cor}
\begin{proof}
  Choose a morphism $\alpha'' \in H_{\alpha'}$ which is not contained
  in any other component of $\Hom_{\alpha'}(X,Z)_{\red}$. Now apply
  Proposition~\ref{prop:p11} to $f'' = \alpha'' \circ \beta$ and any
  component $H_{f''} \subset \Hom_f(X,Y)_{\red}$ that contains $f''$.
\end{proof}

\begin{cor}\label{cor:surjprop}
  The morphism $\eta$ is surjective and proper.
\end{cor}
\begin{proof}
  Since $\Hom_f(X, Y)$ is connected, surjectivity of $\eta$ follows
  from Proposition~\ref{prop:p11}.  Since $\Hom_{\alpha}(X, Z)$ is
  connected, properness of $\eta$ follows from
  Corollary~\ref{cor:prop}.
\end{proof}

\subsection{The max.~étale factorization of a deformed morphism}

Let $f' \in \Hom_f(X,Y)_{\red}$ be any deformation of $f$. The
surjectivity of $\eta$ implies that $f'$ factor via $Z$. Here we will
show that $f'$ has $Z$ as maximally étale factorization. To this end,
let
\begin{equation}
  \label{eq:fact_of_fprime}
  \xymatrix{ {X} \ar[r]_{\alpha'} \ar@/^0.3cm/[rr]^{f'} & {Z'}
    \ar[r]_{\beta'} & {Y}}
\end{equation}
be the maximally étale factorization of $f'$. The universal property
from Definition~\ref{def:maxetale} then yields a morphism $Z' \to Z$.
Reversing the roles of $f$ and $f'$, we also obtain a morphism $Z \to
Z'$ which shows that $Z$ and $Z'$ are isomorphic. \qed

\subsection{Étalité of $\tilde \eta$}

Since surjective and generically injective finite morphisms between
normal spaces are biholomorphic, the following lemma suffices to prove
that for each pair of points in the normalizations, $\tilde \alpha'
\in \widetilde \Hom_{\alpha}(X,Z)$ and $\tilde f' \in \widetilde
\Hom_{f}(X,Y)$ with $\tilde \eta (\tilde \alpha') = \tilde f'$, the
morphism $\tilde \eta$ induces an isomorphism of analytic
neighborhoods. This shows that $\tilde \eta$ is étale and ends the
proof of Theorem~\ref{thm:main0a}.

\begin{lem}
  Let $\alpha'$ be a point in $\Hom_\alpha(X,Z)$ and $f' := \beta
  \circ \alpha'$. Then there are open neighborhoods $U = U(\alpha')$
  and $V = V(f')$ such that $\eta|_U : U \to V$ is bijective.
\end{lem}
\begin{proof}
  Let $y \in Y$ be a general point, and $\Omega = \Omega(y)$ a
  sufficiently small analytic neighborhood such that 
  $$
  \beta^{-1}(\Omega) = \Omega_{1,Z} \cup \cdots  \cup \Omega_{n,Z} \text{\quad and \quad }
  (f')^{-1}(\Omega) = \Omega_{1,X} \cup \cdots  \cup \Omega_{n\cdot m,X}
  $$
  are disjoint unions of open sets which are each isomorphic to
  $\Omega$. If $\Omega' \subset \subset \Omega$ is a relatively
  compact neighborhood of $y$, the sets
  \begin{align*}
    (U_i)_{1\leq i \leq n} & := \{ \alpha'' \in \Hom_{\alpha}(X, Z)_{\red} \,|\, \alpha''(\Omega'_{1,X}) \subset \Omega_{i,Z} \}  \\
    V & := \{ f'' \in \Hom_f(X, Y)_{\red} \,|\, f''(\Omega'_{1,X})
    \subset \Omega \}
  \end{align*}
  are open, the $U_i$ are disjoint, and $\eta^{-1}(V) = \cup_{1 \leq i
    \leq n} U_i$. Using that $\beta|_{\Omega_{i,Z}} : \Omega_{i,Z} \to
  \Omega$ are biholomorphic, the identity principle then immediately
  implies that $\eta|_{U_i}: U_i \to V$ is injective.
  Proposition~\ref{prop:p11} implies that for any given number $i$,
  $U_i$ is either empty or surjects onto $V$. The proof is finished if
  choose $i$ such that $\alpha' \in U_i$ and set $U := U_i$.
\end{proof}

\section{Infinitesimal decomposition of the Hom-scheme}

Theorem~\ref{thm:main} asserts that a cover of $\Hom(X,Y)$ decomposes
into a torus and deformations that are vertical with respect to the
rational quotient. In this section we will show an infinitesimal
version of the decomposition. We believe that this is of independent
interest.

Before we formulate the result in Theorem~\ref{infini} below, recall
the following standard fact of algebraic group theory.

\begin{fact}\label{fact:maxabelian}
  Let $G$ be an algebraic group. Then there exists a maximal compact
  Abelian subgroup, i.e., an Abelian variety $T \subset G$ which is a
  subgroup and such that no intermediate subgroup $T \subset S \subset
  G$, $T \not = S$, is an Abelian
  variety.
  
  A maximal compact Abelian subgroup is unique up to conjugation.
\end{fact}

The decomposition result then goes as follows.

\begin{thm}\label{infini}
  Let $f: X \to Y$ be a surjective morphism between normal
  complex-projective varieties, and
  $$
  \xymatrix{ {X} \ar[r]_{\alpha}
    \ar@/^0.3cm/[rr]^{f} & {Z} \ar[r]_{\beta} & {Y}} 
  $$
  be the maximally étale factorization of $f$. Then there is a
  canonical decomposition of the space of infinitesimal deformations
  of $f$,
  $$
  T_{\Hom(X,Z)}|_f = \Hom \left(f^*(\Omega^1_Y),\O_X \right) =
  \mathfrak a \oplus V,
  $$
  where $\mathfrak a \subset H^0(Z,T_Z)$ is the Lie algebra of a
  maximal compact Abelian variety $T \subset \Aut^0(Z)$ and where $V
  \subset \Hom \left(f^*(\Omega^1_Y),\O_X \right)$ is a subspace of
  the space of infinitesimal deformations that are vertical with
  respect to the maximal rationally connected fibration of $Z$.
\end{thm}

Recall that ``infinitesimal deformations that are vertical with
respect to the maximal rationally connected fibration'' were defined
in Definition~\ref{def:vertical} on page~\pageref{def:vertical}.

\begin{rem}
  The functoriality of the maximal rationally chain connected
  fibration, \cite[thm.~ IV.5.5]{K96}, implies that an infinitesimal
  deformation $\sigma \in \Hom \left(f^*(\Omega^1_Y),\O_X \right)$
  that is vertical with respect to the maximal rationally connected
  fibration of $Z$ is also vertical with respect to the maximal
  rationally connected fibration of $Y$.
\end{rem}

\begin{cor}\label{cor:infini}
  In the setup of Theorem~\ref{infini}, if $g \in \Hom(X,Z)_{\red}$,
  then the tangent space $T_{\Hom(X,Z)_{\red}}|_g$ is spanned by
  infinitesimal deformations that vertical with respect to the maximal
  rationally connected fibration of $Z$, and by tangent vectors of the
  $T$-orbit through $g$. \qed
\end{cor}

We prove Theorem~\ref{infini} in the remainder of the present section.
As usual, we subdivide the proof into steps.

\subsection{Reduction to the case of a finite morphism}

Using Stein factorization of the morphism $f$, we can assume without
loss of generality that $f$ ---and hence $\alpha$--- are actually
finite.  In fact, if $f$ is not finite, consider the Stein
factorization as in Diagram~\eqref{diag:stein} on
page~\pageref{diag:stein}: $f = h \circ g$, where $g: X \to W_0$ has
connected fibers and $h: W_0 \to Y$ is finite. For the reduction, we
need to show that the canonical pull-back morphism
$$
g^*: {\Hom}(h^*(\Omega^1_Y),\O_{W_0}) \to
{\Hom}(f^*(\Omega^1_Y),\O_X).
$$
is bijective. Since $g$ is surjective, injectivity is obvious.
Concerning surjectivity of this map, consider an element $u \in
{\Hom}(f^*(\Omega^1_Y),\O_X)$,
$$
u: f^*(\Omega^1_Y) = g^* h^*(\Omega^1_Y) \to \O_X.
$$
The composition of the canonical map $h^*(\Omega^1_Y) \to
g_*g^*h^*(\Omega^1_Y)$ and the push-forward of $u$,
$$
h^*(\Omega^1_Y) \xrightarrow{} g_*g^*h^*(\Omega^1_Y) \xrightarrow{g_*(u)} g_* \O_X  = \O_{W_0}
$$
then yields a morphism $v: h^*(\Omega^1_Y) \to \O_{W_0}$ such that
$g^*(v)$ and $u$ agree over the smooth part of $Y$, where the
pull-back of $\Omega^1_Y$ is locally free. Since the $\Hom$-sheaves
are torsion free, this implies that $u = g^*(v)$.

In summary, we have shown that $a$ is an isomorphism. The reduction
step is then clear.

\subsection{Setup and Notation}

For convenience, let $X_{\Sing}, Y_{\Sing}$ and $Z_{\Sing}$ denote the
singular loci, and set
\begin{align*}
  X^0 & := X \setminus (X_{\Sing} \cup \alpha^{-1}(Z_{\Sing}) \cup f^{-1}(Y_{\Sing})) \quad \text{and} \quad \\
  Z^0 & := Z \setminus (Z_{\Sing} \cup \alpha(X_{\Sing}) \cup
  \beta^{-1}(Y_{\Sing})).
\end{align*}
Then $\codim_X (X\setminus X^0) = \codim_Z (Z\setminus Z^0) \geq 2$.
The space of infinitesimal deformations can thus be rewritten as
follows.
\begin{equation}
  \label{eq:decom}
  \begin{aligned}
    T_{\Hom(X,Y)}|_f & = \Hom(f^*(\Omega^1_Y),\O_X) \\
    & = \Hom(f^*(\Omega^1_Y)^{\vee\vee},\O_X) = \Hom(f^*(\Omega^1_Y) \vert X^0,\O_{X^0}) \\
    & = H^0(X^0,f^*(T_Y)) = H^0(Z^0, \alpha^*(T_Z)) \\
    & = H^0(Z^0, \alpha_*(\O_{X^0}) \otimes T_{Z^0}) = \Hom_{Z^0}(\alpha_*(\O_{X^0})^*,T_{Z^0}). \\
  \end{aligned}
\end{equation}

\begin{notation}
  If $\sigma \in H^0(X, f^*(T_Y))$ is any infinitesimal deformation of
  the morphism $\alpha$, let $\hat \sigma \in
  \Hom_{Z^0}(\alpha_*(\O_{X^0}), T_{Z^0})$ be the associated morphism.
\end{notation}

\begin{lem}\label{lem:wrtg}
  Let $\sigma$ be an infinitesimal deformation. Then
  \begin{equation*}
  \sigma|_{X^0} \in H^0\left(X^0, \alpha^* \Image(\hat \sigma)\right).    
  \end{equation*}
\end{lem}
\begin{proof}
  The claim immediately follows from the definition of $\hat \sigma$:
  if $z \in Z^0$ is a general point, and $\alpha^{-1}(z) = \{x_i
  |i=1\ldots m\}$, then the image of $\hat \sigma$ at $z$ is spanned
  by the tangent vectors $T\alpha(\sigma(x_i))_{i=1\ldots m}$.
\end{proof}

\subsection{Decomposition of the Infinitesimal Deformations}

Recall Fact~\ref{fact:splitting} which asserts that
$\alpha_*(\O_{X^0}) \cong \O_{Z^0} \oplus \sE_{\alpha}^\vee$.  This,
together with the Equations~\eqref{eq:decom} yields a decomposition
\begin{equation}
  \label{eq:decomp}
  \Hom(f^*(\Omega^1_Y),\O_X)  = H^0(Z, T_Z) \oplus \underbrace{\Hom_{Z^0}(\sE_{\alpha}, T_{Z^0})}_{=: V'}.
\end{equation}

\begin{notation}
  If $\sigma \in H^0(X, f^*(T_Y))$ is any infinitesimal deformation of
  the morphism $\alpha$, let $\hat \sigma' \in H^0(Z, T_Z)$ and $\hat
  \sigma'' \in \Hom_{Z^0}(\sE_{\alpha}, T_{Z^0})$ be the associated
  vector field and morphism, respectively.
\end{notation}

\subsection{Interpretation of $V'$}

We will now show that infinitesimal deformations $\sigma$ of $\alpha$,
which correspond to elements in $V'$ are vertical with respect to the
rational quotient of $Z$. To this end, choose an ample bundle $H \in
\Pic(Z)$ and let $C \subset Z$ be an associated general complete
intersection curve. Fact~\ref{fact:splitting} and the characterization
of the maximally étale factorization, Corollary~\ref{cor:EAntiAmple},
then assert that the restriction $\sE_{\alpha}^{\vee}|_C$ is
anti-ample. It follows that $\sE_{\alpha}|_C$ is ample, and so is its
image in $T_{Z^0}|_C$ under the map $\hat \sigma''$. The refinement of
Miyaoka's characterization of uniruled manifolds, Fact~\ref{fact:HN3},
implies that $\Image(\hat \sigma'')$ is then vertical with respect to
the rational quotient of $Z$, and Lemma~\ref{lem:wrtg} yields the
claim.

\subsection{The Abelian variety $T$  and end of the proof of Theorem~\ref{infini} } 

We consider the connected algebraic group $\Aut^0(Z)$. By a classical
theorem of Chevalley, there exists an extension
$$
0 \to L \to \Aut^0(Z) \to T' \to 0
$$
where $L$ is linear-algebraic and $T'$ an Abelian variety. This
sequence is not necessarily split, but it is known
\cite[thm.~3.12]{Lib78} that there is a maximal compact Abelian
subgroup $T \subset \Aut^0(Z)$ such that the induced map $T \to T'$ is
étale.

Let $\mathfrak a \subset H^0(Z,T_Z)$ be the subalgebra generated by
$T$ and $\mathfrak a'$ that one generated by $L$. This gives a
decomposition
$$
H^0(Z,T_Z) = \mathfrak a \oplus \mathfrak a'.
$$
Since $L$ is linear-algebraic, the closures of its orbits are
rationally connected. As a consequence, $L$ acts trivially on the
rational quotient $Q_Z$, hence $\mathfrak a'$ is vertical and we
obtain a decomposition
$$
\Hom \left(f^*(\Omega^1_Y),\O_X \right) = \mathfrak a \oplus V
$$
with $V := \mathfrak a' \oplus V'$ vertical. This ends the proof of
Theorem~\ref{infini}.  \qed

\section{Decomposition of the Hom-scheme, Proof of Theorem~\ref{thm:main}}
\label{sec:decomp}

The proof of Theorem~\ref{thm:main}, which we give in this section, is
the longest and most involved in this paper. Before we start with all
the details in Section~\ref{sec:startofproof} below, we give a short
idea of proof.

\subsection{Idea of proof}
\label{sec:decomp-idea}

In Section~\ref{sec:startofproof} we will quickly reduce to the case
where $f$ is finite. For simplicity, assume further that the maximally
rationally connected fibration $q_Y: Y \dasharrow Q_Y$ is a morphism
and that the maximally étale factorization is an isomorphism.  Let $T
\subset \Aut^0(Y)$ be a maximal compact Abelian subgroup, as in
Fact~\ref{fact:maxabelian} above.

Under these assumptions, the composition morphism
$$
\begin{array}{rccc}
\tau: & \Hom_f(X,Y) & \to & \Hom(X, Q_Y) \\
& f' & \mapsto & q_Y \circ f'
\end{array}
$$
is equivariant with respect to the natural $T$-action on
$\Hom_f(X,Y)$ and $\Hom(X, Q_Y)$, respectively. The infinitesimal
decomposition of the $\Hom$-scheme, Theorem~\ref{infini}, then asserts
that the image of $\tau$ contains a dense $T$-orbit. By properness,
the image of $\tau$ will be homogeneous under the $T$-action.

The standard fact that actions of Abelian varieties on rationally
connected varieties are necessarily trivial (note that this is not true
for rationally chain connected varieties!)
then implies that
$T$-orbits in $\Hom_f(X,Y)$ surject finitely onto the image of $\tau$,
better still, that they are étale over the image of $\tau$. This
quickly gives the decomposition.

The main difficulty in the proof of Theorem~\ref{thm:main} is that
$q_Y$ need not be regular. Although the space of rational maps $X
\dasharrow Q_Y$ can easily be defined as a subscheme of $\Hilb(X
\times Q_Y)$, its universal properties are too weak to construct a
morphism similar to $\tau$ above ---see \cite{Hanamura87, Hanamura88}
for a discussion of the complications that already arise with the
space of birational automorphisms. We will need to consider a somewhat
weaker construction instead.

\subsection{Reduction to the case of a finite morphism}
\label{sec:startofproof}

Using the stability of Stein factorization under deformation,
Proposition~\ref{prop:stability-of-stein}, we can assume without loss
of generality that the morphism $f$ is finite. Throughout, we consider
the maximally étale factorization of $f$,
$$
\xymatrix{ X \ar[r]_{\alpha} \ar@/^0.3cm/[rr]^{f} & Z
  \ar[r]_{\beta} & Y, }
$$
where $\beta$ is étale in codimension 1.

\subsection{Setup of notation}

Before we seriously start the proof of Theorem~\ref{thm:main}, we need
to set up some notation.

\begin{notation}
  Let $q_Z : Z \dasharrow Q_Z$ and $q_Y : Y \dasharrow Q_Y$ be the
  maximal rationally connected fibrations and
  $$
  \nu : \widetilde \Hom_f(X,Y) \to \Hom_f(X,Y)_{\red}
  $$
  be the normalization morphism.  Let $T \subset \Aut^0(Z)$ be a
  maximal compact Abelian subgroup, as discussed in
  Fact~\ref{fact:maxabelian} above.
\end{notation}

\begin{rem}\label{rem:equivariance1}
  The group $T \subset \Aut^0(Z)$ naturally acts on $Z$ and on
  $\Hom_{\alpha}(X,Z)_{\red}$. By Proposition~\ref{prop:canonrc} and
  the conventions fixed in Notation~\ref{not:alwayscanaon}, the group
  $T$ acts also on the base of the maximally rationally connected
  fibration $Q_Z$. With these actions, the maximally rationally
  connected fibration map $q_Z : Z \dasharrow Q_Z$ is automatically
  $T$-equivariant wherever it is defined.
\end{rem}

As a next step, we define subvarieties $H_{\vertical}^g \subset
\Hom_{\alpha}(X,Z)$ which are the analogues to the fibers of the map
$\tau$ that was discussed in the introductory
Section~\ref{sec:decomp-idea} above.

\begin{notation}
  If $g \in \Hom_{\alpha}(X,Z)$ is any morphism, define the reduced
  subvariety
  $$
  H_{\vertical}^g := \{ h \in \Hom_{\alpha}(X,Z)_{\red} \,|\, q_Z \circ
  g = q_Z \circ h \}.
  $$
  As in Notation~\ref{not:reldef1}, we call $H_{\vertical}^g$ the
  ``space of relative deformations of $g$ over $q_Z$''. Consider the
  restricted group action morphism
  $$
  \begin{array}{rccc}
    \mu_g : & T \times H^g_{\vertical} & \to & \Hom_{\alpha}(X,Z)_{\red} \\
    &(t,\alpha') & \mapsto &  t \circ \alpha'
  \end{array}
  $$
\end{notation}

\begin{rem}\label{rem:equivariance2}
  If $t \in T$ and $g \in \Hom_{\alpha}(X,Z)_{\red}$ are any two
  elements, the associated vertical deformation spaces of $g$ and $t
  \cdot g$ differ only by translation in $\Hom_{\alpha}(X,Z)_{\red}$.
  More precisely, we have $H_{\vertical}^{t \cdot g} = t \cdot
  H_{\vertical}^g$. This follows trivially from the equivariance of
  $q_Z$.
\end{rem}

\subsection{Study of the restricted action morphism}

The spaces $\Hom_{\alpha}(X,Z)_{\red}$ and $H^{\alpha}_{\vertical}$
are not necessarily proper. We will show now, however, that the
restricted group action map $\mu_{\alpha}$ is still a proper morphism.
This will suffice to prove both parts of Theorem~\ref{thm:main}.

\begin{prop}\label{prop:72}
  The restricted action morphism $\mu_{\alpha} : T \times
  H^{\alpha}_{\vertical} \to \Hom_{\alpha}(X,Z)_{\red}$ is proper and
  surjective. It becomes étale after passing to the normalization.
\end{prop}

Assume for the moment that Proposition~\ref{prop:72} holds true. We
will first show that this implies Theorem~\ref{thm:main} and then, in
Sections~\ref{sec:pp721}--\ref{sec:pp722} below, prove the
proposition.

\begin{proof}[Proof of Theorem~\ref{thm:main}, Statement (\ref{thm:main}.\ref{thm:m1})]
  Let $\tilde \mu_{\alpha} : T \times \tilde H^{\alpha}_{\vertical}
  \to \widetilde \Hom_{\alpha}(X,Z)$ be the étale morphism between the
  normalizations that is associated with $\mu_{\alpha}$.  Let 
  $$
  \eta : \Hom_{\alpha}(X, Z)_{\red} \to \Hom_{f}(X, Y)_{\red}
  $$
  be the proper and surjective composition morphism discussed in
  Theorem~\ref{thm:main0a}, and $\tilde \eta$ the associated étale
  morphism between the normalizations. By Proposition~\ref{prop:72}
  and Theorem~\ref{thm:main0a}, the composition
  $$
  \hat \mu_{\alpha} := \tilde \eta \circ \tilde \mu_{\alpha} : T \times \tilde
  H^{\alpha}_{\vertical} \to \widetilde \Hom_f(X,Y)
  $$
  is then surjective and étale, and it suffices to show that
  \begin{equation}\label{eq:rctorc}
    (\eta \circ \mu_{\alpha}) (\{e\} \times H^{\alpha}_{\vertical}) = 
    \eta(H^{\alpha}_{\vertical}) \subset H^f_{\vertical}.
  \end{equation}
  For this, observe that the universal property of the maximally
  rationally chain connected fibration, \cite[IV~thm.~5.5]{K96}, shows
  the existence of a commutative diagram of dominant rational maps as
  follows.
  \begin{equation}
    \label{eq:univprop}
    \xymatrix{
      X \ar[r]_{\alpha} \ar@/^.3cm/[rr]^{f} & Z \ar[r]_{\beta} \ar@{-->}[d]^{q_Z} & Y \ar@{-->}[d]^{q_Y}\\
      & Q_Z \ar@{-->}[r]_{\beta_Q} & Q_Y}
  \end{equation}
  Equation~\eqref{eq:rctorc} then follows by definition of $H$ and
  $H^{\alpha}_{\vertical}$.
\end{proof}

\begin{proof}[Proof of Theorem~\ref{thm:main}, Statement (\ref{thm:main}.\ref{thm:m2})]
  If $f$ is maximally étale, i.e., if $\beta$ is isomorphic, the claim
  follows trivially from the construction.
  
  Now assume that $Y$ is smooth. We will see below that
  $H^{\alpha}_{\vertical}$ is a connected component of
  $\eta^{-1}(H^f_{\vertical})$. It is, however, generally false that
  $\eta^{-1}(H^f_{\vertical}) = H^{\alpha}_{\vertical}$. The map $\hat \mu_{\alpha}
  = \tilde \eta \circ \tilde \mu_{\alpha}$ does therefore not always
  satisfy statement (\ref{thm:main}.\ref{thm:m2}) of
  Theorem~\ref{thm:main} and needs to be modified accordingly, by
  adding more components to $T \times \tilde H^{\alpha}_{\vertical}$,
  one for each connected component of $\eta^{-1}(H^f_{\vertical})$.
  More precisely, we will show that there exists a finite set $T_R
  \subset T$, such that $\eta^{-1}(H^f_{\vertical}) \subset
  \Hom_{\alpha}(X, Z)$ is the disjoint union of copies of
  $H^{\alpha}_{\vertical}$ that are realized in $\Hom_{\alpha}(Z,Y)$
  as translates of $H^{\alpha}_{\vertical}$ by elements of $T_R$ under
  the natural $T$-action on $\Hom_{\alpha}(Z,Y)$, i.e.
  \begin{equation}
    \label{eq:decoPH}
    \eta^{-1}(H^{\alpha}_{\vertical}) = \overset{\bullet}{\bigcup_{t \in T_R}} t \cdot H^{\alpha}_{\vertical}.
  \end{equation}
  We can therefore consider the modified restricted action morphism
  $$
  \begin{array}{rccc}
    \mu_{\alpha}' : & T \times ( T_R \times H^{\alpha}_{\vertical} ) & \to & \Hom_{\alpha}(X,Z)_{\red} \\
    & \bigl(t_1, (t_2, \alpha ')\bigr) & \mapsto & t_1 \circ t_2 \circ \alpha'
  \end{array}
  $$
  It is then obvious that the associated morphism $\tilde
  \mu_{\alpha}'$ between normalizations is étale.  Setting $\mu :=
  \tilde \eta \circ \tilde \mu_{\alpha}'$ then finishes the proof.
  
  It remains to find $T_R$. To this end, we need to introduce the
  following two subgroups of $T$.
  \begin{align*}
    T_{\vertical, Z} & := \{ t \in T \,|\, q_Z = q_Z \circ t \} \\
    T_{\vertical, Y} & := \{ t \in T \,|\, \beta_Q \circ q_Z = \beta_Q \circ q_Z \circ t \}
  \end{align*}

  \begin{claim}\label{claim:Tvert}
    The subgroups $T_{\vertical, Z}$ and $T_{\vertical, Y}$ are both finite.
  \end{claim}

  \begin{proof}[Proof of the claim]
    Since $Y$ is smooth, the quotient map $q_Y$ is almost holomorphic
    in the sense discussed in Fact~\ref{fact:mrcf}.  The general
    $q_Y$-fiber $Y_q \subset Y$ is thus smooth, rationally connected
    and therefore \cite[cor.~4.18]{Debarre01} simply connected. Recall
    that $\beta$ is étale in codimension 1, i.e.~étale away from a set
    of codimension $\geq 2$.  Zariski-Nagata's theorem on the purity
    of the branch locus, \cite[thm.~3.1]{SGA1} implies that $\beta$ is
    étale.  The preimage $\beta^{-1}(Y_q)$ is then a disjoint union of
    several copies of the rationally connected manifold $Y_q$, each a
    fiber of $q_Z$.  This observation has two consequences.
    
    First, the well-known fact that actions of connected,
    positive-dimensional Abelian varieties on rationally connected
    manifolds must necessarily be trivial, \cite[lem.~5.2]{Fujiki78},
    implies that $T_{\vertical, Z}$ is discrete, hence finite.
    
    Second, the observation shows that the dominant rational map
    $\beta_Q$ defined in Diagram~\eqref{eq:univprop} is generically
    finite. This implies that $T_{\vertical, Y}$ is finite.
    Claim~\ref{claim:Tvert} is thus shown.
  \end{proof}
  
  To apply Claim~\ref{claim:Tvert}, recall from
  Theorem~\ref{thm:main0a} and Proposition~\ref{prop:72} that any
  morphism $f' \in \Hom_f(X,Y)_{\red}$ can be decomposed as $f' =
  \beta \circ t \circ \alpha'$ where $\alpha' \in
  H^{\alpha}_{\vertical}$ and $t \in T$. We then have equivalences
  \begin{align*}
    &&  f' & \in H^f_{\vertical} \\
    \Leftrightarrow &&    q_Y \circ f' & = q_Y \circ f && \text{Definition} \\
    \Leftrightarrow &&    q_Y \circ \beta \circ t \circ \alpha' & = q_Y \circ \beta \circ \alpha && \text{Diagram~\ref{eq:univprop}} \\
    \Leftrightarrow &&    \beta_Q \circ q_Z \circ t \circ \alpha' & = \beta_Q \circ q_Z \circ \alpha' && \text{since $\alpha' \in H^{\alpha}_{\vertical}$} \\
    \Leftrightarrow && \beta_Q \circ q_Z \circ t & = \beta_Q \circ q_Z && \text{because $\alpha'$ is surjective} \\
    \Leftrightarrow && t & \in  T_{\vertical,Y}  && \text{Definition}
  \end{align*}
  This already shows that
  $$
  \eta^{-1}(H^f_{\vertical}) = \bigcup_{t \in T_{\vertical, Y}} t \cdot H^{\alpha}_{\vertical}
  $$
  It follows immediately from the definition that two translates,
  $t_1 \cdot H^{\alpha}_{\vertical}$ and $t_2 \cdot
  H^{\alpha}_{\vertical}$ are equal if and only if $t_1\cdot t_2^{-1} \in
  T_{\vertical, Z}$, and otherwise disjoint. We can therefore take
  $T_R$ to be a system of representatives for the finite group
  quotient $\factor T_{\vertical, Y} . T_{\vertical, Z}.$.
  
  Assuming that Proposition~\ref{prop:72} holds, this ends the proof
  of Theorem~\ref{thm:main}.
\end{proof}

\subsection{Proof of Proposition~\ref{prop:72}, a rational decomposition of the $\Hom$-scheme}
\label{sec:pp721-epsilon}

The aim of this section is to construct a rational analogue of the
function $\tau'$ defined in the introduction. To simplify the
notation, we consider the irreducible components of
$\Hom_{\alpha}(X,Z)_{\red}$ separately.

\begin{notation}
  Let
  $$
  \Hom_{\alpha}(X,Z)_{\red} = \bigcup_i \Hom_{\alpha}(X,Z)_{\red}^i
  $$
  be the decomposition into irreducible components, and let
  $$
  H^g_{\vertical,i} := H^g_{\vertical} \cap \Hom_{\alpha}(X,Z)_{\red}^i
  $$
  be the associated decomposition of the $H^g_{\vertical}$.
\end{notation}

By definition of $H^g_{\vertical, i}$, the space
$\Hom_{\alpha}(X,Z)_{\red}^i$ is naturally decomposed into a disjoint
union of subvarieties,
\begin{equation}
  \label{eq:decompset}
  \Hom_{\alpha}(X,Z)_{\red}^i = \overset{\bullet}{\bigcup_{g \in \Hom_{\alpha}(X,Z)_{\red}^i}} H^g_{\vertical, i}
\end{equation}
where all subvarieties $H^g_{\vertical, i}$ are all fibers of the
set-theoretic map
$$
\begin{array}{rccc}
\tau'_i : & \Hom_{\alpha}(X,Z)_{\red}^i & \to & \{ \text{rational maps } X \dasharrow Q_Z \}. \\
& g & \mapsto & q_Z \circ g
\end{array}
$$
We have already discussed in Section~\ref{sec:decomp-idea} that it
might not be possible to define a good scheme-structure on the set of
rational maps which makes $\tau'_i$ a morphism. To construct an
algebraic substitute for $\tau'_i$, observe that $\mathbb C$ is
uncountable. Equation~\eqref{eq:decompset} therefore decomposes
$\Hom_{\alpha}(X,Z)_{\red}^i$ into uncountable many subvarieties. If
$\overline{\Hom_{\alpha}(X,Z)_{\red}^i}$ is a projective
compactification of $\Hom_{\alpha}(X,Z)_{\red}^i$, then
$\Chow(\overline{\Hom_{\alpha}(X,Z)_{\red}^i})$ has only countably
many components. The decomposition of $\Hom_{\alpha}(X,Z)_{\red}^i$
therefore yields a rational map between varieties
$$
\begin{array}{rccc}
\tau_i : & \Hom_{\alpha}(X,Z)_{\red}^i & \dasharrow & \Chow(\overline{\Hom_{\alpha}(X,Z)_{\red}^i}) \\
& g & \mapsto & \text{ closure of } H^g_{\vertical, i}
\end{array}
$$
that agrees with $\tau'_i$ on an open subset.

Although $\tau_i$ is only a rational map, there is a little that we
can say about its infinitesimal structure.

\begin{lem}\label{lem:vertandtau}
  If $\Hom_{\alpha}(X,Z)_{\red}^i \subset \Hom_{\alpha}(X,Z)_{\red}$
  is any irreducible component and $g \in \Hom_{\alpha}(X,Z)_{\red}^i$
  a general point, then the kernel of the tangent morphism
  $\ker(T\tau_i|_g)$ is exactly the space of vertical infinitesimal
  deformations.
  
  In particular, the tangent space
  $T_{\Hom_{\alpha}(X,Z)_{\red}^i}|_g$ at $g$ is spanned by the
  tangent space to $H^g_{\vertical, i}$, and by the tangent space to
  the $T$-orbit through $g$.
\end{lem}
\begin{proof}
  Define a distribution, i.e.~a saturated subsheaf $\sF \subset
  T_{\Hom_{\alpha}(X,Z)_{\red}^i}$ of $\O_X$-modules, as follows: if
  $h \in \Hom_{\alpha}(X,Z)_{\red}^i$ is a smooth, general point, let
  $$
  \sF|_h \subset T_{\Hom_{\alpha}(X,Z)_{\red}^i}|_h \subset
  T_{\Hom_{\alpha}(X,Z)}|_h = \Hom(h^*(\Omega^1_Z), \O_X)
  $$
  be those elements that are
  \begin{itemize}
  \item tangent to the reduced scheme $\Hom_{\alpha}(X,Z)_{\red}^i$,
    and
  \item vertical with respect to the rational quotient $q_z : Z
    \dasharrow Q_Z$.
  \end{itemize}
  
  It is clear that $\ker(T\tau_i|_g) \subset \sF|_g$. To show the
  other inclusion, assume that we are given a vertical infinitesimal
  deformation $\vec v \in \sF|_g$. In order to prove
  Lemma~\ref{lem:vertandtau}, we need to show that $\vec v$ is tangent
  to $H^g_{\vertical, i}$.  For this, consider a holomorphic arc
  $\gamma: \Delta \to \Hom_{\alpha}(X,Z)_{\red}^i$ such that
  \begin{enumerate}
  \item $\gamma(0)=g$, and the derivatives satisfy
  \item $\gamma'(0)= \vec v$, and 
  \item $\gamma'(t) \subset \gamma^*(\sF)$, for all $t
  \in \Delta$.
  \end{enumerate}
  The infinitesimal description of the universal morphism $\Delta
  \times X \to Z$ then shows that the image of $\gamma$ is entirely
  contained in $H^g_{\vertical, i}$. It follows that $\vec v \in
  T_{H^g_{\vertical, i}}|_g = \ker(T\tau_i|_g)$.
\end{proof}

\subsection{Proof of Proposition~\ref{prop:72}, properness and surjectivity}
\label{sec:pp721}

With the preparations from the previous section, we can now prove the
first assertion of Proposition~\ref{prop:72}. We show surjectivity
first for the restricted action morphism $\mu_{\tilde g}$, where
$\tilde g$ is a general element.

\begin{lem}\label{lem:gensurjectivity}
  If $\Hom_{\alpha}(X,Z)_{\red}^i \subset \Hom_{\alpha}(X,Z)_{\red}$
  is any irreducible component and $\tilde g \in
  \Hom_\alpha(X,Z)_{\red}^i$ a general point, then $\mu_{\tilde g}$
  surjects onto $\Hom_{\alpha}(X,Z)_{\red}^i$.
\end{lem}
\begin{proof}
  Lemma~\ref{lem:vertandtau} and the infinitesimal decomposition,
  Corollary~\ref{cor:infini}, together imply that $\mu_{\tilde g}$ is
  of maximal rank at $\tilde g$ and therefore dominates the component
  $\Hom_{\alpha}(X,Z)_{\red}^i$.
  
  To show that $\mu_{\tilde g}$ is surjects onto
  $\Hom_{\alpha}(X,Z)_{\red}^i$, it suffices to show that its image is
  closed, i.e.~that that the limit of every convergent sequence of
  points in the image is again contained in the image. Using the
  compactness of $T$, this follows immediately.
\end{proof}

It is now easy to extend the surjectivity result to all $g \in
\Hom_\alpha(X,Z)_{\red}^i$.

\begin{lem}\label{lem:surj2}
  In the setup of Lemma~\ref{lem:gensurjectivity}, if $g \in
  \Hom_\alpha(X,Z)_{\red}^i$ is any point, then $\mu_g$ and
  $\mu_{\tilde g}$ have the same image in $\Hom_{\alpha}(X,Z)_{\red}$.
\end{lem}
\begin{proof}
  The surjectivity of $\mu_{\tilde g}$,
  Lemma~\ref{lem:gensurjectivity}, implies that there exist elements
  $t \in T$ and $\hat g \in H_{\vertical,i}^{\tilde g}$ such that
  $t\cdot \hat g=g$. By Remark~\ref{rem:equivariance2}, we have
  $H^g_{\vertical} = t\cdot H_{\vertical}^{\tilde g}$ and therefore
  $$
  \mu_g = t \circ \mu_{\tilde g} \circ (Id, t^{-1}).
  $$
  This shows the claim.
\end{proof}

\begin{cor}\label{cor:surjofmualpha}
  If $g \in \Hom_\alpha(X,Z)_{\red}$ is any point, then $\mu_g$ is
  surjective and proper. In particular, $\mu_\alpha$ is surjective and
  proper.
\end{cor}
\begin{proof}
  The surjectivity of $\mu_g$ follows immediately from
  Lemma~\ref{lem:surj2} and the fact that $\Hom_\alpha(X,Z)_{\red}$ is
  connected by definition.
  
  It remains to show that $\mu_g$ is proper, i.e.~that the preimage of
  any compact set $K \subset \Hom_\alpha(X,Z)_{\red}$ is again
  compact. But again, given a sequence $(t_n, g_n) \subset
  \mu_g^{-1}(K)$, using that $T$ is compact and the sequence $t_n\cdot
  g_n$ has a cumulation point in $K$, it is easy to prove that $(t_n,
  g_n)$ has a convergent subsequence.
\end{proof}

\subsection{Proof of Proposition~\ref{prop:72}, étalité}
\label{sec:pp722}

The étalité of $\tilde \mu_{\alpha}$ will be deduced using the
following criterion. Although fairly standard, we found no reference
in the literature and give a quick proof.

\begin{prop}[Étalité criterion]\label{prop:etalite}
  Let $f : X \to Y$ be a proper, finite morphism between irreducible
  varieties and assume that $Y$ is normal. If there exists a number $d
  \in \mathbb N$ such that for all $y \in Y$, the preimages
  $f^{-1}(y)$ contains (set-theoretically) exactly $d$ points, then
  $f$ is étale.
\end{prop}
\begin{proof}
  Let $y \in Y$ be any point and $f^{-1}(y) = \{x_1, \ldots, x_d\}$.
  By \cite[sect.~2.3]{CAS} we can find an analytic neighborhoods $U$
  of $y$ and $V_i$ of $x_i$ such that $f^{-1}(U) = V_1 \cup \cdots
  \cup V_d$ and such that the $V_i$ are disjoint. The restrictions
  $f|_{V_i} : V_i \to U$ must then be bijective and, by the analytic
  version of Zariski's main theorem, \cite[prop.~14.7]{Remmert94},
  biholomorphic. This shows the claim.
\end{proof}

\begin{proof}[Proof of Proposition~\ref{prop:72}]
  We have already seen in Corollary~\ref{cor:surjofmualpha} that
  $\mu_\alpha$ is proper and surjective. It remains to show that the
  associated morphism between the normalizations is étale.
  
  For this, let $\widetilde \Hom_\alpha(X,Z)_{\red}$ and $\tilde
  H^\alpha_{\vertical}$ be the normalizations of
  $\Hom_\alpha(X,Z)_{\red}$ and $H^\alpha_{\vertical}$, respectively.
  Further, let
  $$
  \tilde \mu_{\alpha} : T \times \tilde H^\alpha_{\vertical} \to
  \widetilde \Hom_\alpha(X,Z)_{\red}
  $$
  be the morphism associated with $\mu_{\alpha}$. This morphism
  will then also be proper.
  
  By the étalité criterion, Proposition~\ref{prop:etalite}, it remains
  to show that the number of elements in fibers of $\tilde
  \mu_{\alpha}$ is constant.
  
  Recall that $T$ acts effectively and freely on
  $\Hom_\alpha(X,Z)_{\red}$, and therefore freely on the normalization
  $\widetilde \Hom_\alpha(X,Z)_{\red}$. If $G \subset T$ denotes the
  ineffectivity of the $T$-action on $Q_Z$, i.e.~the kernel of the
  natural map $T \to \Aut(Q_Z)$, then $G$ acts freely on
  $H^\alpha_{\vertical}$ and $\tilde H^\alpha_{\vertical}$. Here we
  need to consider the natural $G$-action on $T \times \tilde
  H^\alpha_{\vertical}$, where $G$ acts on the factor $T$ by left
  multiplication.  This action is likewise free.
  
  Proposition~\ref{prop:72} is shown if we prove that for any pair
  $(t,\tilde g) \in T \times \tilde H_{\vertical}^\alpha$, the
  associated $\tilde \mu_\alpha$-fiber is exactly the $G$-orbit, i.e.
  $$
  \tilde \mu_\alpha^{-1} \bigl( \tilde \mu_\alpha (t,\tilde g) \bigr) = G \cdot (t,\tilde g).
  $$
  The inclusion ``$\supseteq$'' is clear.
  
  For the other inclusion, consider two pairs contained in the same
  fiber,
  \begin{equation}\label{eq:samefiber}
    \tilde \mu_\alpha (t_1,\tilde g_1) = \tilde \mu_\alpha (t_2,\tilde g_2)  
  \end{equation}
  If $\nu : \widetilde \Hom_\alpha(X,Z)_{\red} \to
  \Hom_\alpha(X,Z)_{\red}$ is the normalization morphism,
  equation~\eqref{eq:samefiber} then implies
  $$
  \nu(\tilde g_1) = t_1^{-1}t_2 \cdot \nu(\tilde g_2) = \nu(t_1^{-1}t_2 \cdot \tilde g_2)
  $$
  The assumption $\tilde g_1, \tilde g_2 \in \tilde
  H_{\vertical}^\alpha$, i.e.~ $q_Z \circ \nu(\tilde g_1) = q_Z \circ
  \nu(\tilde g_1) = q_Z \circ \alpha$ then yields that $t_1^{-1}t_2
  \in G$, which ends the proof of Proposition~\ref{prop:72} and hence
  of Theorem~\ref{thm:main}.
\end{proof}

\providecommand{\bysame}{\leavevmode\hbox to3em{\hrulefill}\thinspace}
\providecommand{\MR}{\relax\ifhmode\unskip\space\fi MR}
\providecommand{\MRhref}[2]{%
  \href{http://www.ams.org/mathscinet-getitem?mr=#1}{#2}
}
\providecommand{\href}[2]{#2}

\end{document}